\documentclass[letterpaper, 10 pt, conference]{ieeeconf}  % Comment this line out
\usepackage{color}
\usepackage{textcomp}
\usepackage[mathcal]{euscript}
\usepackage{mathrsfs}
\usepackage{amsfonts,amssymb,amsbsy}
\usepackage{graphicx,empheq}
\usepackage{nomencl}
\usepackage{multirow}
\usepackage[tableposition=top]{caption}
\def\e1{{\varepsilon_{1}}}
\def\b1{{\beta_{1}}}
\def\bp3{{\beta_{3}}}
\def\ep3{{\varepsilon_{3}}}

\newcommand{\mb}{\mathbf}

\MHInternalSyntaxOn
\providecommand*\phantomword[3][c]{%
\mathchoice
{\MT_phantom_word:NNnn #1\displaystyle {#2}{#3}}%
{\MT_phantom_word:NNnn #1\textstyle {#2}{#3}}%
{\MT_phantom_word:NNnn #1\scriptstyle {#2}{#3}}%
{\MT_phantom_word:NNnn #1\scriptscriptstyle {#2}{#3}}%
}
\def\MT_phantom_word:NNnn #1#2#3#4{%
\@begin@tempboxa\hbox{$\m@th#2#4$}%
\setlength\@tempdima{\widthof{$\m@th#2#3$}}%
\hbox{\hb@xt@\@tempdima{\csname bm@#1\endcsname}}%
\@end@tempboxa}
\MHInternalSyntaxOff

\IEEEoverridecommandlockouts                              % This command is only
\title{\LARGE \bf
Modeling an elastic beam with piezoelectric patches by including magnetic effects
}

\author{A. \"Ozkan \"Ozer and K. A. Morris % <-this % stops a space
\thanks{ This  research was supported by a Discovery Grant from the Natural Sciences and Engineering Research Council of Canada (NSERC) and by U.S.  AFOSR grant FA9550-10-1-0530.}% <-this % stops a space
\thanks{A. \"{O}zkan \"{O}zer is a postdoctoral researcher in the Department of Mathematics \& Statistics, University of Nevada, Reno, NV, 89523,  USA
        {\tt\small aozer@unr.edu}}%
        \thanks{K. Morris is a professor in the Department of Applied Mathematics, University of Waterloo, ON, N2L3G1,  Canada
        {\tt\small kmorris@uwaterloo.ca}}%
}

\begin{document}

\maketitle
\thispagestyle{empty}
\pagestyle{empty}

%%%%%%%%%%%%%%%%%%%%%%%%%%%%%%%%%%%%%%%%%%%%%%%%%%%%%%%%%%%%%%%%%%%%%%%%%%%%%%%%

\begin{abstract}
%\chg{Rewritten. note removed ref. to wave eqn- EB model has a beam equation}
Models for piezoelectric beams using Euler-Bernoulli small displacement theory predict the dynamics of slender beams at the low frequency accurately but  are insufficient for beams vibrating at high frequencies or beams with low length-to-width aspect ratios. A more thorough model that includes the effects of rotational inertia and  shear strain,  Mindlin-Timoshenko small displacement theory,  is needed to predict the dynamics  more accurately for these cases. Moreover, existing models ignore the magnetic effects  since the magnetic effects are relatively small. However, it was shown recently  \cite{O-M1} that these effects can substantially change the controllability and stabilizability properties of even a single piezoelectric beam. In this paper, we use a variational approach to derive  models that include magnetic effects  for an elastic beam with two piezoelectric patches actuated by different voltage sources.  Both Euler-Bernoulli and Mindlin-Timoshenko small displacement theories are considered.  Due to the magnetic effects, the equations are quite different from the standard equations.
 %and therefore voltage controls arise from the jump conditions at the edges of the patches.

\end{abstract}

\section{Introduction}
\begin{figure}[h!tb]
\centering
\includegraphics[width=2.7in]{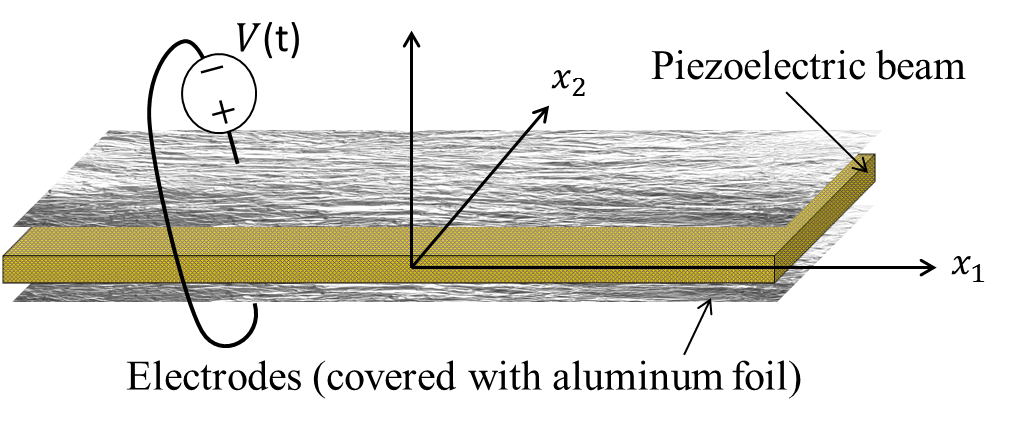}
\caption{\footnotesize When voltage $V(t)$ is supplied to the electrodes, an electric field is created  and the beam either shrinks or extends.}
\label{piezosingle}
\end{figure}

Piezoelectric materials have been successfully used to transfer mechanical energy into electro-magnetic energy, and vice versa. Due to their low cost, light weight, and high durability, they have been very competitive for many tasks in many fields ranging from space technology to the automotive industry.
%\kchg{ The simplest form of a piezoelectric structure is a single piezoelectric beam. }
A single piezoelectric beam is an elastic beam covered by electrodes on its top and bottom surfaces, insulated at the edges (to prevent fringing effects), and connected to an external  electric circuit. (See Figure \ref{piezosingle}.)

\begin{figure}[htb]
\centering
\includegraphics[width=2.7in]{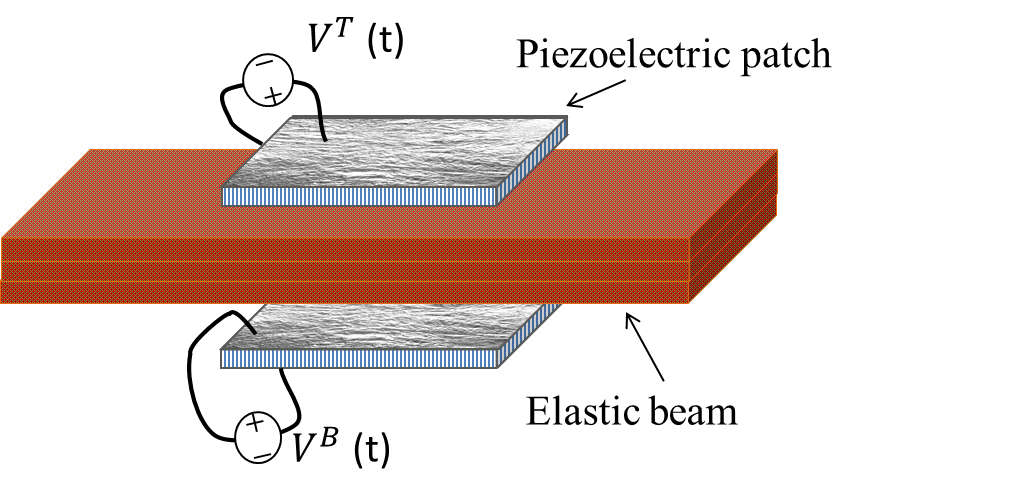}
\caption{\footnotesize An elastic beam with piezoelectric patches actuated by voltage $V^T(t)$ at the top and $V^B(t)$ at the bottom. %Depending on  the sign and the magnitude of $V^T(t)$ and $V^B(t)$, the elastic beam stretches or bends.
 }
\label{fig:patch}
\end{figure}

We use a  variational approach to describe the dynamics for a
piezoelectric beam with both Euler-Bernoulli and Mindlin-Timoshenko small displacement assumptions. Using the same methodology, we describe the dynamics for the common configuration of an elastic beam with two piezoelectric patches bonded on the elastic beam top and bottom surfaces as shown in Figure \ref{fig:patch}.
Mechanical effects are modelled through Euler-Bernoulli (E-B) or Mindlin-Timoshenko (M-T) small displacement assumptions, for instance see \cite{Banks-Smith}, \cite{E-Inman}, \cite{Hansen}, \cite{Smith}, \cite{Yang}.  The (M-T) model includes all the effects of the (E-B) model, and additionally, it includes the effects of the shear strain and rotational inertia. The (E-B) beam model is not sufficient to describe the dynamics of thin beams vibrating at  high  frequencies or of beams with low length-to-width aspect ratios (short and stubby beams), see \cite{Dietl} and the references therein.
%\kchg{moved discussion  of energy space to conclusions. and deleted some since this paper only derives model, not state-space}
%It is also observed in \cite{Hansen} (charge actuation) that the control operator for the piezoelectric model derived by the (M-T) model is smoother in the natural energy space  than the one derived by the (E-B) beam model.
%Here the natural energy space is is chosen so that the states of the model  have enough degree of smoothness in the weak formulation.
Our (M-T) model is different from the model  in \cite{Zhang} where the degree of smoothness of the control operator is the same as that of the (E-B) model.

Most studies of piezoelectric structures consider only mechanical and electrical effects but not magnetic effects.
Electrical and magnetic effects are modeled by Maxwell's equations with three widely used assumptions: electro-static, quasi-static and fully dynamic \cite{Tiersten}. Electro-static and quasi-static approaches are the most widely used, see for instance  \cite{K-M-M2}, \cite{L-M}, \cite{Rogacheva}, \cite{Smith}, \cite{Tiersten}, \cite{Tzou}. Assuming stationary electrical effects is generally considered a reasonable assumption since it has been experimentally observed that  magnetic effects are a very minor aspect of the overall dynamics for polarized ceramics,
 %\kchg{and therefore electrostatic or quasi-static models predict the dynamical behavior; }
 see the review article \cite{Yang1}. However, in \cite{O-M}, it is shown that the controlled piezoelectric beam with magnetic effects is not strongly stabilizable for many  parameter values. This is quite different from the stabilizability and controllability of a piezoelectric beam without magnetic effects.

% The choice of the feedback controller, total current at the electrodes, is very practical and can be measured easily \cite{Ronkanen}.  In \cite{O-M1}, it is proved that a single piezoelectric beam  can never be exactly observable (controllable) in the natural energy space.
%
In this paper, we include all magnetic effects and model a beam-patch system driven by two different voltage sources.  The dynamics of electro-magnetic effects are included. Hamilton's Principle leads to strongly  coupled equations for stretching and bending with the (E-B) assumptions, and for stretching, bending and rotation with the (M-T) assumptions. Conditions at the boundary of the patches are also obtained. %The voltage control at both patches acts through one of the boundary conditions.

%%%%%%%%%%%%%%%%%%%%%%%%%%%%%%%%%%%%%%%%%%%%%%%%%%%%%%%%%%%%%%%%%%%%%%%%%%%%%%%%
\section{A single piezoelectric beam}
In this section we present  the initial and boundary value problem for models of a piezoelectric beam  that include  full magnetic effects.  Throughout this paper, we use dots and primes to denote differentiation with respect to time  $ (t)$ and space $(x),$ respectively.  We first start with the modeling of piezoelectric patches. Let $x_1, x_2$ be the longitudinal directions and let $x_3$ be transverse directions. We assume that the elastic beam  occupies the region $\Omega=[0,L]\times [-r,r]\times [-\frac{h}{2}, \frac{h}{2}].$ We denote by $\partial \Omega$ the boundary of $\Omega;$ the electroded region and the insulated region.

%\chg{This would be clearer and tie in better to the next section if done on [a,b] not [0,L]!!! See the patch section! }

 We adopt the following linear constitutive relationship \cite{Tiersten} for piezoelectric beams
\begin{eqnarray}
\label{cons-eqq10}
\left( \begin{array}{l}
 T \\
 D \\
 \end{array} \right)=
\left[ {\begin{array}{*{20}c}
   c & -\Gamma^{\text{T}}  \\
   \Gamma & \varepsilon  \\
\end{array}} \right]\left( \begin{array}{l}
 S \\
 E \\
 \end{array} \right)
\end{eqnarray}
where $T=(T_{11}, T_{22}, T_{33}, T_{23}, T_{13}, T_{12})^{\text T}$ is the stress vector, $S=(S_{11}, S_{22}, S_{33}, S_{23}, S_{13}, S_{12})^{\text T}$  is the strain vector, $D=(D_1, D_2, D_3)^{\text T}$ and  $E=(E_1, E_2, E_3)^{\text{T}}$ are the electric displacement  and the electric field vectors, respectively, and moreover, the matrices $[c], [\Gamma], [\varepsilon]$ are the matrices with elastic, electro-mechanic and dielectric constant entries. (For more details, see \cite{Tiersten}.) Under the assumption of transverse isotropy and polarization in $x_3-$direction, these matrices reduce to $\varepsilon={\text{diag}}(\varepsilon_{1}, \varepsilon_{2},\varepsilon_{3}),$ and
\begin{eqnarray}
\nonumber &c= \left[ {\begin{array}{*{20}c}
   {c_{11} } & {c_{12} } & {c_{13} } & 0 & 0 & 0  \\
   {c_{21} } & {c_{22} } & {c_{23} } & 0 & 0 & 0  \\
   {c_{31} } & {c_{32} } & {c_{33} } & 0 & 0 & 0  \\
   0 & 0 & 0 &    c_{44} & 0 & 0 \\
      0 & 0 & 0 &    0 & c_{55} & 0 \\
         0 & 0 & 0 &    0 & 0 & c_{66}
\end{array}} \right],\\
\nonumber &\Gamma= \left[ {\begin{array}{*{20}c}
   0 & 0 & 0 &    0 & \gamma_{15} & 0 \\
      0 & 0 & 0 &    -\gamma_{15} & 0 & 0 \\
         \gamma_{31} & \gamma_{31} & \gamma_{33} &    0 & 0 & 0
\end{array}} \right].&
%\nonumber  &\varepsilon= \left[ {\begin{array}{*{20}c}
 %  \varepsilon_{1} & 0 & 0 \\
  %    0 & \varepsilon_{2} & 0\\
   %      0 &  0 & \varepsilon_{3} \\
%\end{array}} \right].&
\end{eqnarray}

In standard beam theory, all forces acting in the $x_2$ direction, and the transverse normal stress $T_{33}$ are negligible. In (E-B) beam theory, the shear strain $ S_{13}$ is zero but it is  nonzero in  (M-T) beam theory.

 Denote flexural rigidity, modulus of elasticity, piezoelectric coefficients, and  dielectric constants by $c_{11}, c_{55}, e_{31}, e_{15}, \e1$ and $\ep3$, respectively.
 Let $v=v(x),$ $w=w(x)$ and $\psi=\psi(x)$ denote the longitudinal displacement of the center line, transverse    displacement , and the rotation of the beam, respectively. For simplicity of notation we set $x=x_1$ and $z=x_3.$ Continuing with the small-displacement assumptions, the displacement fields, strains, and the linear constitutive equations \cite{Tiersten} for the (E-B) and (M-T) beam models are given in Table \ref{const} with the following notation
%\chg{line up on separate lines to be easier to read if there is room}
\begin{eqnarray}
\nonumber \gamma_1= \gamma_{15},~ \gamma_3=\gamma_{31},~\alpha_1=\alpha_{11} + \gamma_3^2 \beta_3,~ \beta_{1}=\frac{1}{\varepsilon_{1}}\\
\label{coef} \alpha_3=\alpha_{11} + \gamma_1^2 \beta_1, ~\alpha_{11}=c_{11},~~\alpha_{33}=c_{55}, ~\beta_3=\frac{1}{\varepsilon_{3}}.
%\nonumber && \beta=\beta_{33}=\frac{1}{\ep3}, ~~ e=e_{31}, ~~ \gamma= e \beta, \\
\end{eqnarray}
\renewcommand{\arraystretch}{1.4}
\begin{table}[h]
\centering
\small
\begin{tabular}{|l|l|l|} \hline
& Euler-Bernoulli  (E-B) \\ \hline
\multirow{2}{*}{Displacement fields} &$U_1=v(x)-zw'(x)$   \\
&$U_3=w(x)$  \\ \hline
\multirow{2}{*}{Strains}&$S_{11}=\frac{\partial U_1}{\partial x} =v'-zw''$    \\
&$S_{13}=\frac{1}{2}\left(\frac{\partial U_1}{\partial z}+ \frac{\partial U_3}{\partial x}\right)=0$\\ \hline
\multirow{4}{*}{Constitutive equations}&$T_{11}=\alpha_1 S_{11}-\gamma_3\beta _3D_3$   \\
&$T_{13}=-\gamma_1 \beta_1 D_1$  \\
&$E_1=\beta_{1}D_1$ \\
&$E_3=-\gamma_3\beta_3 S_{11}+\beta_3 D_3$  \\ \hline
\end{tabular}\\
\vspace{0.2in}
\small
\begin{tabular}{|l|l|l|}
\hline
& Mindlin-Timoshenko (M-T)\\ \hline
\multirow{2}{*}{Displacement fields} & $U_1=v(x)+z\psi(x)$  \\
& $U_3=w(x)$  \\ \hline
\multirow{2}{*}{Strains} & $S_{11}=\frac{\partial U_1}{\partial x}=v'+z\psi'$   \\
& $S_{13}=\frac{1}{2}\left(\frac{\partial U_1}{\partial z}+ \frac{\partial U_3}{\partial x}\right)=\frac{1}{2} (w' + \psi)$\\ \hline
\multirow{4}{*}{Constitutive equations}& $T_{11}= \alpha_1 S_{11}-\gamma_3\beta_3 D_3$   \\
& $T_{13}= \alpha_3\ S_{13}-\gamma_1 \beta_1 D_1$  \\
& $E_1=-\gamma_1 \beta_1 S_{13}+\b1 D_1$ \\
& $E_3=-\gamma_3\beta_3 S_{11}+\bp3 D_3$  \\ \hline
\end{tabular}
\footnotesize
\caption{\footnotesize Displacement fields, strains, and constitutive equations for  two different beam models.
%In the above (E-B) and (M-T) refer to Euler-Bernoulli and Mindlin-Timoshenko.
}
\label{const}
\end{table}

Electro-magnetic effects are described by Maxwell's equations:
  \begin{subequations}
  \label{Maxwell}
\begin{empheq}[left={\phantomword[r]{0}{ }  \empheqlbrace}]{align}
\label{Gauss-law} \nabla\cdot D =~\sigma_b & \quad{\rm{in}} \quad\Omega \times \mathbb{R}^+~& \\
\label{Gauss-magne} \nabla\cdot B=~0 &\quad{\rm{in}} \quad\Omega \times \mathbb{R}^+~ & \\
\label{Faraday} \nabla\times E=~-\dot B & \quad{\rm{in}} \quad\Omega \times \mathbb{R}^+~&\\
\label{Ampere}  \frac{1}{\mu}(\nabla\times B)= ~i_b +  \dot D & \quad{\rm{in}} \quad\Omega \times \mathbb{R}^+~&
\end{empheq}
\end{subequations}
where the dots denote  differentiation with respect to time $t.$
Here $B$ denotes magnetic field vector, and $\sigma_b, i_b, \sigma_s, i_s, V, \mu, n$ denote body charge density, body current density, surface charge density,
 surface current density, voltage, magnetic permeability, and unit normal vector respectively.
In this paper it is assumed that the only external force acting on the beam is the voltage applied at the electrodes and so the essential electric boundary conditions are
\begin{subequations}
\begin{empheq}[left={\phantomword[r]{0}{ }  \empheqlbrace}]{align}
\label{charge_boun}  - D\cdot n =~\sigma_s  & \quad{\rm{on}} \quad {\partial\Omega} \times \mathbb{R}^+~\quad&\text{(Charge )}\\
\label{current_boun} \frac{1}{\mu}(B \times n) =~i_s & \quad{\rm{on}} \quad {\partial\Omega} \times \mathbb{R}^+~\quad & \text{(Current)}\\
\label{voltage_boun} \phi=~V  & \quad{\rm{on}} \quad {\partial\Omega} \times \mathbb{R}^+~\quad& \text{(Voltage).}
\end{empheq}
\end{subequations}
%There are also the appropriate mechanical boundary condition at the edges of the beam  (clamped, hinged, free, etc.).

There are mainly three approaches to inclusion of electromagnetic effects in piezo-electric beams \cite{Tiersten}:
\vspace{0.1in}

\noindent \textbf{(i)} Electrostatic: In this widely-used approach,
  magnetic effects are completely ignored: $B=\dot D=i_b=\sigma_b=0.$  Maxwell's equations  (\ref{Maxwell}) reduce to
   $\nabla\cdot D =0$ and $\nabla \times E=0.$ Therefore, by Poincar\'{e}'s theorem there exists a scalar electric potential
   such that $E=-\nabla \phi$  where $\phi$ is determined up to a constant .\\

\noindent \textbf{(ii)} Quasi-static: This approach rules out some but not all  the magnetic effects:% (non-magnetizable materials) and
  $\sigma_b=i_b=0$, but $\dot D$ and $B$ are non-zero and (\ref{Maxwell}) becomes
       \begin{eqnarray}\nonumber\nabla \cdot D=0, ~\nabla\cdot B=0, ~\dot B= -\nabla\times E, ~ \dot D= \frac{1}{\mu}(\nabla\times B).\end{eqnarray} The equation $\nabla\cdot B=0$
       implies that there exists a magnetic potential vector $A$ such that $B=\nabla\times A.$ It follows from substituting $B$ to $\dot B= -\nabla\times E$ that there exists a scalar electric potential $\phi$ such that $E=-\nabla\phi-\dot A.$ he maTgnetic potential $A$ is not unique (See \cite{O-M1}). %One simplification in this approach is to ignore $A$ and $\dot A$ since $A, \dot A\ll \phi.$  Note that we still have $\dot D$ non-zero.
       \\

\noindent \textbf{(iii)} Fully dynamic: In this approach,  $A$ and $\dot A$ are left in the model.  Depending on the type of material, body charge density $\sigma_b$ and body current density $i_b$ can also be  non-zero.
%Note that even though the  displacement current $\dot D$ is assumed to be non-zero in both quasi-static and fully dynamic approaches, the term $\ddot D$  is  zero in quasi-static approach since $\dot A=0.$
%The electrostatic assumption has been widely used in the literature to reduce the complexity of Maxwell's equations due the minor magnetic effects on the dynamics of piezoelectric structures. %However, these effects become significant when it comes to the controller design.
\vspace{0.1in}

In this paper, the third, fully dynamic approach for modeling piezoelectric beams that includes all of the magnetic effects, is used as in \cite{O-M}.
The magnetic field $B$ is perpendicular to the $x-y$ plane due to Gauss's law of magnetism, i.e. $\nabla\cdot B=~0,$  and therefore
$B=(B_1, B_2, B_3)$ has only the $y-$component $B_2$ nonzero, and it is only a function of $x.$  This is  because the surface current at the electrodes have only $x-$component (tangential)  and $B$ is perpendicular to both the
outward normal vector ($n=(0,0,1)$ or $n=(0,0,-1)$) at the electrodes and the current on the electrodes. For simplicity, we also assume that $E_1=0,$ and thus $D_1=0$ by Table \ref{const}. Maxwell's equations including the effects of $B$ become
$$ \nabla\cdot B=0, ~~~\dot B= -\nabla\times E, ~~~ \dot D= \frac{1}{\mu}(\nabla\times B).$$

From the last equation it follows that  $$\frac{ d B_2}{dx}=-\mu \dot D_3,$$ and so $$B_2=-\mu\int \dot D_3(\xi, x_3, t) ~d\xi.$$
The next assumption is that $D_3$ is constant in $x_3:$ $D_3(x,x_3,t)=D_3(x,t).$
Now let \begin{eqnarray}\label{defp}q=\int  D_3(x,  t) ~dx\end{eqnarray} so that $q'=D_3$ and $-\mu \dot q=B_2.$ The magnetic energy (which can be regarded as the electric kinetic energy) is
\begin{eqnarray}\nonumber \mb B &=&  \frac{1}{2\mu}\int_{\Omega} (B_2)^2~dX =  \frac{\mu}{2}\int_{\Omega} \dot q^2 ~dX.
\end{eqnarray}

Since the voltage is prescribed at the electrodes, we use the  Lagrangian \cite{Lee,O-M}
\begin{eqnarray}\label{tildeL}  \mb{L}= \int_0^T \left[\mb{K}-(\mb{P}+\mb{E})+\mb B +\mb{W}\right]~dt\end{eqnarray}
 where $\mb K,$ $\mb P+\mb E,$ $\mb B,$ and $\mb W$ denote the (mechanical) kinetic energy, total stored energy, magnetic energy (electrical kinetic energy) of the beam,  and the work done by the external forces, respectively. This is  different than that in \cite{Hansen} since that paper considers charge actuation. Using (\ref{defp}) and Table \ref{const}, leads to the energies
{\small{
  \begin{align}
 \nonumber  &\mb P+\mb E=\frac{1}{2}\int_\Omega \left(T_{11}S_{11} + D_3 E_3\right) ~dX\\
\nonumber  &~~=\frac{h}{2} \int_0^L \left[ \alpha_1(  v'^2 + \frac{h^2}{12} w''^2) -2 \gamma_3 \beta_3  v' q' + \beta_3  q'^2\right]dx, \quad \text{(E-B)} \\
 \nonumber & \mb K=\frac{\rho}{2} \int_\Omega \left(\dot U_1^2+ \dot U_3^2\right)~dX\\
\nonumber & ~~= \frac{\rho h}{2} \int_0^L \left(\dot v^2+\frac{ h^2}{12} \dot w'^2 + \dot w^2\right)~dx, \quad \text{(E-B)}\\
%\nonumber  &\frac{1}{2}\int_\Omega \left(T_{11}S_{11} + T_{13}S_{13}+ D_3 E_3\right) ~dX&\\
%\nonumber &\frac{h}{2} \int_0^L \left[ \alpha_1\left(  v'^2 + \frac{h^2}{12} \psi'^2\right)+ \right.&
% &\left.+ \alpha_3 (w'+\psi)^2-2 \gamma_1\beta_1  v' q' + \beta_{3} q'^2\right]~dx,&\text{(M-T)}
 \nonumber  &\mb P+\mb E=\frac{1}{2}\int_\Omega \left(T_{11}S_{11} + T_{13}S_{13}+ D_3 E_3\right) ~dX\\
\nonumber  &~~=\frac{h}{2} \int_0^L \left[ \alpha_1\left(  v'^2 + \frac{h^2}{12} \psi'^2\right)+ \alpha_3 (w'+\psi)^2\right., \\
\nonumber & \quad ~~\left.-2 \gamma_1\beta_1  v' q' + \beta_{3} q'^2-2 \gamma_1\beta_1  v' q' + \beta_{3} q'^2\right]~dx, \quad \text{(M-T)}\\
 \nonumber& \mb K=\frac{\rho}{2} \int_\Omega \left(\dot U_1^2+ \dot U_3^2\right)~dX\\
 \nonumber & ~~= \frac{\rho h}{2} \int_0^L \left(\dot v^2+\frac{ h^2}{12} \dot \psi^2 + \dot w^2\right)~dx, \quad \text{(M-T)}\\
\nonumber &\mb{B}= \frac{1}{2\mu}\int_\Omega \|B\|^2 ~dX= \frac{\mu h}{2}\int_{0}^L \dot q^2~dx,\text{(E-B), \text(M-T)}\quad\\
\label{work-done}& \mb{W}=- \int_{\partial\Omega} D_3 ~\bar\phi ~d\Gamma=-\int_0^L  q'~ V(t)~dx, \text{(E-B), \text(M-T)}
\end{align}
}}
where $\bar \phi$ is the electric potential,  and $V(t)$ is the voltage applied at the electrodes. Note that voltage $V(t)$ is the potential difference between the top and bottom electrodes.% \chg{In the above we skipped the integration in $x_2$ variable since all functions are independent of $x_2.$  }

Application of Hamilton's principle, setting the variation of admissible displacements $\{v,w,q\}$ for (E-B) and $\{v,\psi, w, q\}$ of (M-T) of $\bf L$ to zero, yields three partial differential equations. If we assume that  the beam is free at both ends $x=0,L,$ we obtain the sets of partial differential equations and  boundary conditions for the (E-B) assumptions
 \begin{subequations}
  \label{PDEs1}
\begin{empheq}[left={\phantomword[l]{}{{}}\quad\quad \empheqlbrace }]{align}
  &\rho  \ddot v-\alpha_1   v''+\gamma_3\beta      q'' = 0&  \\
\label{ref1}  &\mu \ddot q-\beta_3   v''+\gamma_3\beta_3      q'' = 0 \quad \text{(Stretching)}&\\
     &\rho h \ddot w +\frac{\alpha_1 h^3}{12}w''''=0 \quad\quad \text{(Bending)} &
      \end{empheq}
 \begin{empheq}[left={\phantomword[c]{}{\small\text{(E-B)~}}\quad\quad \empheqlbrace}]{align}
     & \left.~\alpha_1  v' -\gamma_3\beta_3  q' ~\right|_{x=0,L}=0 & \\
      &\left.~\beta_3  q'  -\gamma_3\beta_3 v' ~\right|_{x=0,L}= -\frac{V(t)}{h}&\\
            &w''(0)=w'''(0)=w''(L)=w'''(L)=0,&
      \end{empheq}
  \end{subequations}
and for the (M-T) assumptions
\begin{subequations}
  \label{PDEs2}
\begin{empheq}[left={\phantomword[l]{}{{}}\quad\quad \empheqlbrace }]{align}
  &\rho  \ddot v-\alpha_1   v''+\gamma_3\beta      q'' = 0&  \\
\label{ref2}  &\mu \ddot q-\beta_3   v''+\gamma_3\beta_3      v'' = 0 \quad \text{(Stretching)} &\\
     &\rho h \ddot w -\alpha_3 h (\psi + w')'=0 \quad\quad \text{(Bending)}&\\
     &\frac{\rho h^3}{12}\ddot \psi  -\frac{\alpha_1 h^3}{12} \psi''&\\
     &\quad\quad\quad + \alpha_3 h (\psi+w')=0\quad \text{(Rotation)}&
      \end{empheq}
 \begin{empheq}[left={\phantomword[c]{}{\small\text{(M-T)~~}} \quad\quad\empheqlbrace}]{align}
     & \left.~\alpha_1  v' -\gamma_3\beta_3  q' ~\right|_{x=0,L}=0 & \\
      &\left.~\beta_3  q'  -\gamma_3\beta_3 v' ~\right|_{x=0,L}= -\frac{V(t)}{h}&\\
            &\left.\psi',~(\psi+w')~\right|_{x=0,L}= 0.&
      \end{empheq}
  \end{subequations}

Since the only external force acting on the beam is the voltage at the electrodes,   the bending and rotation equations are completely decoupled from the stretching equations.
%\begin{rmk}[Another formulation]
%However, by using (\ref{coef}), the above boundary conditions can be simplified to
%$$u_x(0)=u_x(L)=\frac{e_{31}}{c_{11}}\frac{V(t)}{h},~~~~ p_x(0)=p_x(L)=\left(\ep3+ %\frac{e_{31}^2}{c_{11}}\right)\frac{V(t)}{h},~~~\left|w_{xx}=w_{xxx}\right|_{x=0,L}=0$$
%\end{rmk}

%\chg{ Note that our approach is quite different from the quasi-static approach as mentioned earlier. Even though the  displacement current $\dot q$ is assumed to be non-zero in both approaches, the term $\ddot q$ in (\ref{ref1}) and (\ref{ref2})  is assumed to be zero in the quasi-static approach. }

\section{Beam-patch system }
In this section, we consider an elastic beam of length $L$ and height $2h_0$ with two  piezoelectric patches with height $h_1$ bonded one at the top and one at the bottom of the beam. Defining $\omega=[a,b]$ where  $0<a<b<L,$ the symmetrically placed patches occupy the region $\Omega_\omega=\omega\times [-r,r]\times \left([-h_0-h_1, -h_0]\bigcup [h_0, h_0+h_1]\right)\subset \Omega.$  The patches are insulated at the edges, and no external mechanical stress is applied through the edges. They are also assumed to be bonded perfectly so that no slip occurs. Moreover, each patch is covered with electrodes at lower and upper faces (See Figure \ref{fig:patch}). The prescribed voltages at  each patch can be different.
\begin{figure}[h]
\centering
\includegraphics[width=0.50\textwidth]{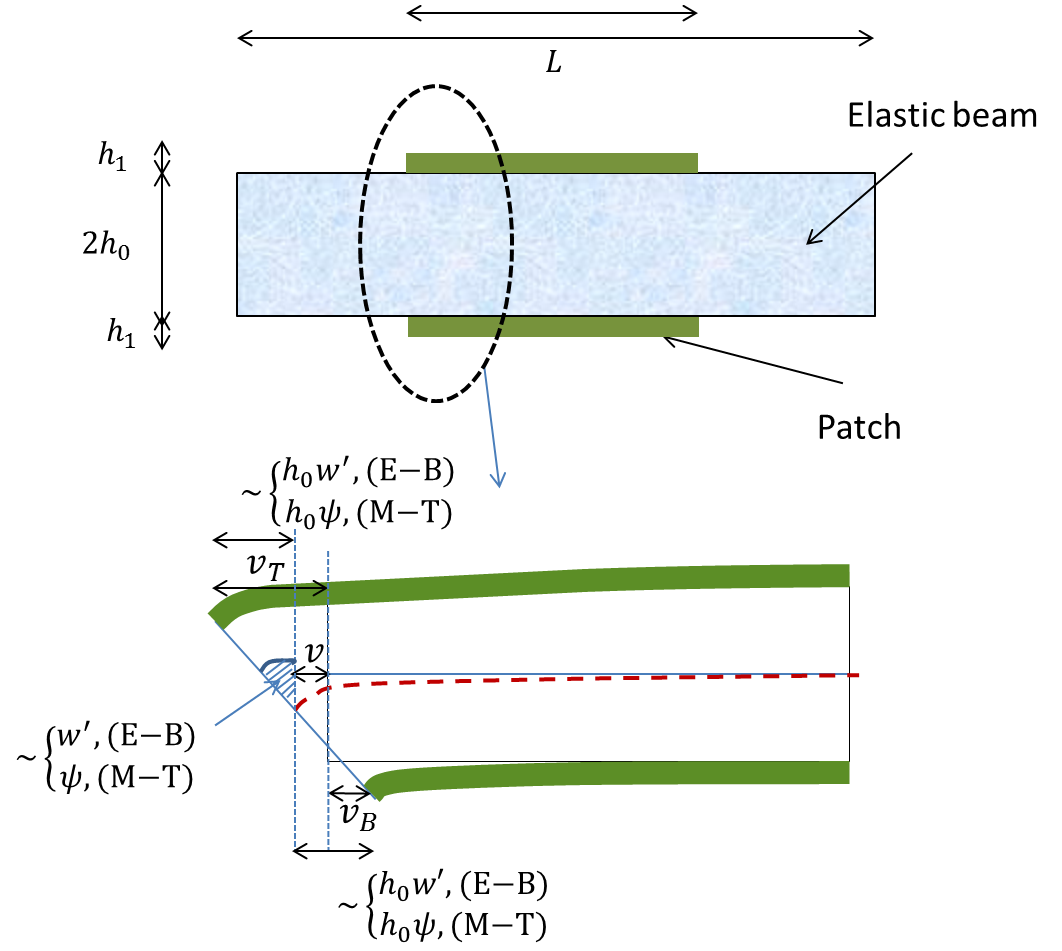}
\caption{\footnotesize Beam-patch system before and after deformation. Due to the geometric constraints we have $\frac{v^T+v^B}{2}=v,~ \frac{v^T-v^B}{2}=h_0w'$  and $\frac{v^T+v^B}{2}=v,~ \frac{v^T-v^B}{2}=h_0\psi$ for the (E-B) and (M-T) assumptions, respectively.}
\label{patch}
\end{figure}

As shown in the previous section (see (\ref{PDEs1}) and (\ref{PDEs2})), when the voltage is prescribed at a patch it either shrinks or extends, so we only consider the stretching motions for the patches.  We use $v^T$ and $v^B$ to denote longitudinal displacements of the centerlines of the top and bottom patches, respectively, and  $v,$ $w,$ and $\psi$ denote the longitudinal displacement of centerline, transverse displacement, and rotation respectively of the elastic beam.  Let $V^T(t)$ and $V^B(t)$ indicate the applied voltages to the top and bottom piezoelectric patches.

% Note that $v^T $ and $v^B$ depends strongly on the prescribed voltages at the electrodes of each patch (see Table \ref{PDEs1} and \ref{PDEs2}).
 Assuming that the patches are bonded firmly and there is no slipping,  we have the following continuity relationship due to the geometry shown in Figure \ref{patch}:
\begin{eqnarray}& %\frac{v^T+v^B}{2}=v, \frac{v^T-v^B}{2}=h_0 w' \Leftrightarrow
v^T=v+h_0 w' ,  v^B = v-h_0w', & \text{(E-B)}\\
& %\frac{v^T+v^B}{2}=v, ~~ \frac{v^T-v^B}{2}=h_0\psi  ~~\Leftrightarrow~~
v^T=v+h_0\psi , ~~ v^B = v-h_0\psi.& \text{(M-T)}
\end{eqnarray}
Letting the indices $b, T, B$ refer to the quantities of  beam, top patch, and bottom patch, respectively,
the total kinetic, potential, electric, magnetic energies and the work done by the external forces
are
\begin{eqnarray}\nonumber \tilde{\mb K}= \mb  K^b + \mb K^T + \mb K^B, ~\tilde{\mb P}= \mb  P^b + \mb P^T + \mb P^B,\\
\nonumber ~\tilde{ \mb E}= \mb E^T + \mb E^B, ~ \tilde{\mb B}= \mb B^T+\mb B^B, ~~ \tilde{\mb W}=\mb W^T + \mb W^B.
\end{eqnarray}
Assume that the material properties of the two patches are identical,
%i.e.,
%\begin{eqnarray}\nonumber && \rho^p=\rho^T=\rho^B, ~\mu=\mu^T=\mu^B, ~\alpha_1^p=\alpha_1^T=\alpha_1^B,\\
%\nonumber &&\alpha_3^p=\alpha_3^T=\alpha_3^B, ~\gamma_1=\gamma_1^T=\gamma_1^B,\\
%\nonumber &&\gamma_3=\gamma_3^T=\gamma_3^B, ~ \beta_1=\beta_1^T=\beta_1^B,~\beta_3=\beta_3^T=\beta_3^B
%\end{eqnarray}
and use the superscript $p$ to indicate the material property of a patch and $b$ to indicate the corresponding property of the beam.
%denotes that the corresponding quantity for the piezoelectric patches.
The magnetic energy and work done are
\begin{align}
\nonumber &\mb B^T= \frac{1}{2\mu}\int_{\Omega_\omega} \|B^T\|^2 ~dX= \frac{\mu h}{2}\int_\omega  (\dot q^T)^2~dx\\
\nonumber &\mb B^B= \frac{1}{2\mu}\int_{\Omega_\omega} \|B^B\|^2 ~dX= \frac{\mu h}{2}\int_\omega(\dot q^B)^2~dx\\
\nonumber &\mb{W}^T= \int_\omega -{q^T}'~ V^T(t)~ dx\\
\nonumber &\mb{W}^B= \int_\omega -{q^B}'~ V^T(t)~ dx. \quad \text{(E-B), (M-T)}\quad\quad
%\nonumber &\gamma_3=\gamma_3^T=\gamma_3^B, ~ \beta_1=\beta_1^T=\beta_1^B,~\beta_3=\beta_3^T=\beta_3^B
\end{align}
%Kinetic energy for the patches are added below. The effect of it was already added to the equations earlier.
The kinetic energy of the beam and patches, potential energy of the beam,  and the total stored energy of patches  are
\begin{align}
\nonumber &\mb K^b=\rho h_0 \int_0^L \left(\dot v^2+\frac{ h^2}{3} \dot w_x^2 + \dot w^2\right)~dx\\
\nonumber &\mb K^T+\mb K^B=\rho h_1 \int_\omega \left[\left(\dot v + h_0 \dot w'\right)^2+\left(\dot v - h_0 \dot w'\right)^2\right]~dx\\
\nonumber &\mb P^b= h_0\int_0^L \left( \alpha_1 \left(  v'^2 + \frac{h_0^2}{3} w''^2\right) \right)~dx\\
\nonumber &\mb P^T + \mb E^T =\frac{h_1}{2} \int_\omega  \left(\alpha_1^p (v'+h_0 w'')^2 \right.\\
\nonumber & \left. \quad\quad\quad\quad\quad\quad -2\gamma_3\beta_3 (v'+h_0 w'') (q^T)' + \beta_3(q^T)'^2\right)~dx \\
\nonumber &\mb P^B + \mb E^B=\frac{h_1}{2} \int_\omega  \left(\alpha_1^p (v'-h_0 w'')^2\right.\\
 \nonumber &\left.-2\gamma_3\beta_3 (v'-h_0 w'') (q^B)' + \beta_3(q^B)'^2\right)~dx, ~~\text{(E-B)}\quad\quad
\end{align}
and with (M-T) assumptions,
\begin{align}
\nonumber &\mb K^b=\rho h_0 \int_0^L \left(\dot v^2+\frac{ h^2}{3} \dot \psi^2 + \dot w^2\right)~dx\\
\nonumber &\mb K^T+\mb K^B=\frac{\rho h_1}{2} \int_\omega\left[\left(\dot v^2+ h_0 \dot \psi\right)^2+\left(\dot v^2- h_0 \dot \psi\right)^2\right]~dx\\
\nonumber &\mb P^b= h_0 \int_0^L \left( \alpha_1\left(  v'^2 + \frac{h_0^2}{3} \psi'^2\right) + \alpha_3 (w'+\psi)^2\right)~dx\\
\nonumber &\mb P^T + \mb E^T =\frac{h_1}{2} \int_\omega \left(\alpha_1^p (v'+h_0\psi')^2 \right.\\
\nonumber &\left.\quad\quad\quad\quad\quad\quad -2\gamma_3\beta_3 (v^T)' (q^T)' + \beta_3(q^T)'^2\right)~dx\\
\nonumber &\mb P^B + \mb E^B=\frac{h_1}{2} \int_\omega  \left(\alpha_1^p (v^B)'^2 \right.\\
\nonumber & \left.\quad\quad\quad\quad\quad-2\gamma_3\beta_3 v^B (q^B)' + \beta_3(q^B)'^2\right)~dx. \quad\text{(M-T)}\quad\quad
\end{align}
The Lagrangian corresponding to the whole beam-patch system (similar ro (\ref{tildeL})) is
\begin{eqnarray}
\label{Lagrangian} \tilde{\mb{L}}=&& \int_0^T \left[\tilde{\mb{K}}-(\tilde{\mb{P}} + \tilde{\mb{E}})+\tilde{\mb{B}}+\tilde{\mb{W}}\right]~dt.
\end{eqnarray}
We set the variation of $\tilde{\mb L}$ with respect to all admissible displacements $\{ v, w, q^T, q^B\}$  for (E-B) and $\{ v, w, \psi, q^T, q^B\}$ for (M-T)  to zero.
Since the bonding is assumed to be perfect, we also use the constraints that $v, w$ and $w'$ for (E-B) and $v,w$ and $\psi$ for (M-T) are continuous at the edges of the patches.
Note that  $q^T$ and $q^B$ are  zero outside the region $\omega.$
Letting $\chi_{\omega}(x)$ indicate the characteristic function of the interval $\omega=(a,b),$ define
$$\rho(x)=h_1\rho^p\chi_{\omega}(x)+h_0\rho,~~  \alpha(x)=h_1 \alpha_1^p\chi_{\omega}(x)+h_0\alpha_1$$
$$\tilde \rho(x)=h_1h_0^2\rho^p\chi_{\omega}(x)+\frac{\rho h_0^3}{3}, ~~A(x)=h_1h_0^2 \alpha_1^p\chi_{\omega}(x)+\frac{\alpha_1 h_0^3}{3}.$$
 Omitting the moment of inertia term $\ddot w''$ in  the (E-B)  model, we obtain the coupled system of partial differential equations
{\small{
 \begin{align}
  \label{BVP1}
\nonumber  &\rho(x) \ddot v-(\alpha(x) v')'+ \frac{\gamma_3\beta_3  h_1}{2} \left(\left( (q^T)' + (q^B)'\right)\chi_{\omega}\right)'   = 0\\
\nonumber  &\tilde\rho(x) \ddot w+ (A(x) w'')''\\
 \nonumber &\quad\quad\quad\quad+ \frac{\gamma_3\beta_3 h_0h_1}{2} \left(\left( (q^T)' - (q^B)'\right)\chi_{\omega}\right)''=0 ~ {\rm in}~ \Omega\\
\nonumber  &  \mu h_1 \ddot q^T -\beta_3 h_1 (q^T)'' + \gamma_3\beta_3 h_1 (v''+h_0 w''')=0 \\
      & \mu h_1 \ddot q^B  - \beta_3 h_1 (q^B)'' + \gamma_3\beta_3 h_1 (v''-h_0w''')=0 ~ {\rm in}~ \omega .
\end{align}}}
%\kchg{  Reworded to remove word "jump" . Also remove the statement of continuity of v etc, this is done above. They are assumptions, not  consequences of the variational method.}
and also  conditions at  the boundary $\partial\omega= a,b $ of the patches where the material is discontinuous. For the (E-B) model they are, %in addition to the previously stated assumption that $v,w$ and $w^\prime$ are continuous,
{\small{
  \begin{align}
\label{BC-nat-patch-3}
\nonumber &\left[\alpha(x)  v'- \frac{\gamma_3\beta_3 h_0 h_1}{2} \left( (q^T)' + (q^B)'\right)\right]_{\partial\omega}=0 ~~  \text{(Lat. force)}\\
\nonumber &\left[ A(x) w''' - \frac{\gamma_3\beta_3 h_0h_1}{2} \left( (q^T)' - (q^B)'\right)'\right]_{\partial\omega}=0~~ \text{(Shear)}\\
\nonumber &\left[ A(x) w'' - \frac{\gamma_3\beta_3 h_0h_1}{2} \left( (q^T)' - (q^B)'\right)\right]_{\partial\omega}=0 ~~ \text{(Ben. moment)}\\
\nonumber &\left[\beta_3 h_1 (q^T)' -\gamma_3\beta_3 h_0h_1 ( v' + h_0  w'')\right]_{\partial\omega}=-V^T ~~  \text{(Volt- Top)}\\
\nonumber &\left[\beta_3 h_1 (q^B)' -\gamma_3\beta_3 h_0h_1 ( v' - h_0  w'')\right]_{\partial\omega}=-V^B~~ . \text{(Volt- Bottom)}
%\\
 %& w' ~{\rm is ~ continuous  } ~ \text{on}~ \partial\omega,  \quad \text{(E-B)}
\end{align}
}}
%\chg{BCs could be all kinds of combinations. just write in one set BCs- c f? }
The boundary conditions at the ends $x=0,L$ can be chosen to be clamped-free
{\small{
  \begin{equation}
\label{BC-nat-patch-1}  v (0)= w (0)= w'(0)= v'(L)= w''(L)= w'''(L)=0.
%\nonumber & v (0)=w(0)= w'(0)= v (L)= w'(L)= w'(L)=0& \text{(c-c)}\\
%\nonumber & v (0)= w (0)= w''(0)= v (L)= w (L)= w''(L)=0 & \text{(h-h)}\\
%\nonumber & v (0,t)= w (0)= w''(0)= v'(L)= w''(L)= w'''(L)=0 & \text{(h-f),}
\end{equation}
}}
For the (M-T) model, the partial differential equations are
{\small{
\begin{align}
  \label{BVP2}
\nonumber  &\rho(x) \ddot v-(\alpha(x) v')'+ \frac{\gamma_3\beta_3  h_1}{2} \left(\left( (q^T)' + (q^B)'\right)\chi_{\omega}\right)'   = 0\\
\nonumber  & \rho \ddot w -\alpha_3 (\psi + w')' =0, \\
\nonumber    & \tilde\rho(x) \ddot \psi -(A(x)\psi')' + \alpha_3 h_0 (\psi+w')\\
 \nonumber &\quad\quad\quad+ \gamma_1\beta_1 h_1\left[\left((q^T)'-(q^B)'\right)\chi_{\omega}\right]' =0 \quad {\rm in}~ \Omega \\
\nonumber    &\mu h_1 \ddot q^T -\beta_3 h_1 (q^T)'' + \gamma_3\beta_3 h_1 (v''+h_0 \psi'')=0 \\
      & \mu h_1 \ddot q^B  - \beta_3 h_1 (q^B)'' + \gamma_3\beta_3 h_1 (v''-h_0\psi'')=0 ~~ {\rm in}~ \omega
\end{align}
}}
with  conditions at the boundaries of the patches
{\small{
  \begin{align}
\label{BC-nat-patch-8}
\nonumber &\left[ \alpha(x) v'- \frac{\gamma_3\beta_3 h_0 h_1}{2} \left( (q^T)' + (q^B)'\right)\right]_{\partial\omega}=0~\text{(Lat. force)}\\
\nonumber &\left[ A(x) \psi' - \frac{\gamma_3\beta_3 h_0h_1}{2} \left( (q^T)' - (q^B)'\right)\right]_{\partial\omega}=0~ \text{(Ben. moment)}\\
\nonumber &\left[\beta_3 h_1 (q^T)' -\gamma_3\beta_3 h_0h_1 ( v' + h_0 \psi')\right]_{\partial\omega}=-V^T(t)~\text{(Volt-Top)}\\
\nonumber &\left[\beta_3 h_1 (q^B)' -\gamma_3\beta_3 h_0h_1 ( v' - h_0  \psi')\right]_{\partial\omega}=-V^B(t)  .~ \text{(Volt-Bottom)}
%\\
%& w'~ {\rm and }~ \psi~ {\rm are ~ continuous  } ~ \text{on}~ \partial\omega. \quad \text{(M-T)}\
\end{align}
}}
There are also boundary conditions at $x=0,L$; for instance
%chg{BCs could be all kinds of combinations. just write in one set BCs - c f? }
clamped-free boundary conditions are
{\footnotesize{
  \begin{equation}
\label{BC-nat-patch-2}
 v (0)= w (0)= \psi(0)= v'(L)= (w'+\psi)(L)= \psi'(L)=0.
%\nonumber & v (0)=w(0)= \psi(0)= v (L)= w(L)= \psi(L)=0 &\text{(c-c)}\\
%\nonumber & v (0)= w (0)= \psi'(0)= v (L)= w (L)= \psi'(L)=0  &\text{(h-h)}\\
%\nonumber & v (0,t)= w (0)= \psi'(0)= v'(L)= (w'+\psi)(L)= \psi'(L)=0  &\text{(h-f).}
\end{equation}
}}
%\kchg{These boundary conditions at the edges of patches  $\partial\omega$ are often referred to as {\em jump conditions} }
Observe that if $V^T=V^B$, there is only purely stretching motion and bending cannot be controlled. Similarly, $V^T=-V^B$ leads to pure bending (and rotation)  motions and stretching is not controlled. In general, all motions - bending, rotation, and stretching - are coupled.

If there are no magnetic effects, so $\mb B_T=\mb B_B=0$, solve the last two equations in (\ref{BVP1}) for $(q^T)''$ and $(q^B)''$ in $\omega$ and substitute them back to the first two equations to obtain
{\small{
  \begin{align}
\nonumber  &\rho(x) \ddot  v -\left(\left(\alpha(x) - \gamma_3^2\beta_3 h_1 \chi_{\omega}\right) v'\right)' = \frac{\gamma_3(V^T+V^B)}{2}\left(\chi_{\omega}\right)'' & \\
\nonumber  & \tilde\rho(x) \ddot  w  + \left(\left(A(x)- \gamma_3^2 \beta_3 h_0h_1  \chi_{\omega}(x)\right)  w''\right)''\\
 \label{pdes-elim-1} &\quad\quad\quad\quad\quad\quad\quad\quad\quad=\frac{h_0\gamma_3 (V^B-V^T)\left(\chi_{\omega}\right)''}{2}, \quad \text{(E-B)}
\end{align}
and
\begin{align}
\nonumber  & \rho \ddot w -\alpha_3 (\psi + w')' =0\\
\nonumber  &\rho(x) \ddot  v -\left(\left(\alpha(x) - \gamma_3^2\beta_3 h_1 \chi_{\omega}\right) v'\right)' = \frac{\gamma_3(V^T+V^B)\left(\chi_{\omega}\right)'}{2} & \\
 \nonumber & \tilde\rho(x) \ddot  \psi  - \left(\left(A(x)- \gamma_3^2 \beta_3 h_0h_1  \chi_{\omega}\right)  \psi'\right)'\\
 \label{pdes-elim-2} & \quad\quad\quad\quad\quad\quad\quad\quad=\frac{h_0\gamma_3 (V^B-V^T)}{2}\left(\chi_{\omega}\right)',  \quad \text{(M-T)}
\end{align}}}
and if we use the relationship $\alpha_1^p=\alpha_1^p + \gamma^2 \beta$ by (\ref{coef}) the boundary value problems obtained above reduce to, in the case of the (E-B) beam model,
\begin{align}
\nonumber  &\rho(x) \ddot  v -\left(\left(h_1 \alpha_{11}^p \chi_{{\omega}}+h_0\alpha_1\right) v'\right)' \\
\nonumber & \quad\quad\quad\quad\quad= \frac{\gamma_3 (V^T+V^B)}{2}\left(\chi_{{\omega}}\right)' \\
\nonumber   & \tilde\rho(x) \ddot  w  + \left(\left(h_1h_0^2 \alpha_{11}^p\chi_{{\omega}}+\frac{\alpha_1 h_0^3}{3}\right)  w''\right)''  \\
 \label{pdes-elim-3} & \quad\quad\quad=-\frac{h_0\gamma_3}{2} (V^T-V^B)\left(\chi_{{\omega}}\right)''\quad \text{(E-B)}
 \end{align}
with boundary conditions at the patch edges
{\small{
  \begin{align}
\nonumber &\left[ \left(h_1 \alpha_{11}^p\chi_{{\omega}}+h_0\alpha_1\right) v'\right]_{\partial\omega}=0~~ \text{(Lat. force)}\\
\nonumber &\left[ \left(h_1h_0^2 \alpha_{11}^p\chi_{{\omega}}+\frac{\alpha_1 h_0^3}{3}\right) w''\right]_{\partial\omega}=0~~\text{(Bending moment)}\\
\nonumber &\left[ \left(\left(h_1h_0^2 \alpha_{11}^p\chi_{{\omega}}+\frac{\alpha_1 h_0^3}{3}\right) w''\right)' \right]_{\partial\omega}=0~~\text{(Shear)}.
%\label{BC-nat-patch-7} &  w'~~ {\rm {is~ continuous ~on}}~~ \partial\omega.
\end{align}
}}
The system of equations (\ref{pdes-elim-3}) coincides with those in \cite{Banks-Smith} and \cite{Smith} obtained using the same physical assumptions and a Newtonian approach. Again, if $V^T=-V^B,$ only the equation for bending is affected by the voltage, and if  $V^T=V^B $ only stretching is affected.
If (M-T) assumptions are used, the boundary value problems reduce to
\begin{align}
\nonumber  & \rho \ddot w -\alpha_3 (\psi + w')' =0 \\
\nonumber  &\rho(x) \ddot  v -\left(\left(h_1 \alpha_{11}^p \chi_{{\omega}}+h_0\alpha_1\right) v'\right)' = \\
\nonumber & \quad\quad\quad\quad\quad\quad\quad\quad\quad \frac{\gamma_3(V^T+V^B)\left(\chi_{{\omega}}\right)'}{2},  \\
   \nonumber & \tilde\rho(x) \ddot  \psi  - \left(\left(h_1h_0^2 \alpha_{11}^p\chi_{{\omega}}+\frac{\alpha_1 h_0^3}{3}\right)  \psi'\right)' \\
   \label{pdes-elim-5}  &\quad\quad\quad\quad\quad\quad=-\frac{h_0\gamma_3(V^T-V^B)\left(\chi_{{\omega}}\right)'}{2} \quad \text{(M-T)}
\end{align}
with the boundary conditions at the patch edges
{\small{
  \begin{align}
\nonumber &\left[ \left(h_1 \alpha_{11}^p\chi_{{\omega}}+h_0\alpha_1\right) v'\right]_{\partial\omega}=0~~ \text{(Lat. force)}\\
\nonumber &\left[ \left(h_1h_0^2 \alpha_{11}^p\chi_{{\omega}}+\frac{\alpha_1 h_0^3}{3}\right) \psi'\right]_{\partial\omega}=0 ~~\text{(Bending moment)} .
%\label{BC-nat-patch-6}   &  w' ~{\text{and}}~ \psi~~ {\rm {are~ continuous ~on}}~~ \partial\omega .
\end{align}
}}

%\chg{Is this done? If not, remove. I  put something in conclusions as future work}
%\kchg{These systems of equations can be shown to have a unique solution using a standard variational technique such as used in \cite{Lagnese-Lions}.}
%As a side note, to reduce the volume of the paper we have not included the well-posedness of these models. The variational approach gives a theoretical base to %show the well-posed of these models  in appropriate Sobolev spaces.

%\chg{Here is how the main theorem goes: Define the spaces $$V_0:=L^2(0,L),\quad V_1:=\left\{z\in H^1(0,L): z(0)=0\right\}$$ $$V_2:=\left\{z\in V_1: %z'(0)=0\right\}$$
%$$\mathcal{W}_1:=V_1 \times  V_2 \times (H^1(\omega))^2 \times V_1^* \times V_2^* \times (L^2(\omega))^2$$
%$$\mathcal{W}_2:=(V_1)^3 \times (H^1(\omega))^2 \times \times (L^2(\omega))^5.$$
%Let $Y_1=(v,w, q^T,q^B,\dot v, \dot w, {\dot q}^T, {\dot q}^T)^{\text T}$ and $Y_2=(v,w,\psi,q^T,q^B, \dot v, \dot w, \dot\psi,  {\dot q}^T, {\dot q}^B)^{\text %T}.$
%\begin{thm} Let $V(t)\in L^2(0,T).$
%
 % \noindent (i) For every $Y_1^0\in \mathcal{W}_1,$  the system  (\ref{BVP1}), (\ref{BC-nat-patch-1}), (\ref{BC-nat-patch-3}) with $Y_1(0)=Y_1^0$ has a unique solution satisfying $Y\in C^0([0,T],\mathcal{W}_1).$

 % \noindent  (ii) For every $Y_2^0\in \mathcal{W}_2,$ the system  (\ref{BVP2}), (\ref{BC-nat-patch-2}), (\ref{BC-nat-patch-8}) with $Y_2(0)=Y_2^0$ has a unique solution satisfying $Y\in C^0([0,T],\mathcal{W}_2).$
%\end{thm}
%}

\section{Conclusion}

%It is also observed in \cite{Hansen} (charge actuation) that the control operator for the piezoelectric model derived by the (M-T) model is smoother in the natural energy space  than the one derived by the (E-B) beam model.
%Here the natural energy space is is chosen so that the states of the model  have enough degree of smoothness in the weak formulation.

%\kchg{Rewrote the conclusions slightly}
In this paper, we have used a variational approach to derive two models, one with the Euler-Bernoulli (E-B) theory and one with the Mindlin-Timoshenko (M-T) theory, for the dynamics in a voltage-controlled beam-patch system. Fully dynamic  magnetic effects are included.  Boundary conditions at the edges of the patches are also obtained.  Magnetic effects account for the wave behavior of the electro-magnetic effects on the patches. In contrast to classical models, the electro-magnetic components of the system are not decoupled from the mechanical components. The beam equations with the (M-T) theory include the shear and rotational effects. This is very important for beams vibrating at the high frequencies where the (E-B) assumptions are not sufficient for predicting the dynamics.

The variational approach used here will facilitate showing well-posedness of the models, using either a port-Hamiltonian approach \cite{O-M} or standard state-space methods \cite{O-M1}. The variational approach leads to a natural definition of the state space in terms of the beam energies, the natural energy space. This is the topic of current work. Note that  the control term in the (E-B) model is less smooth (in the natural energy space) than the one corresponding to the (M-T) model since the rotations in the (E-B) beam model are defined by $\psi=w'$ (see (\ref{pdes-elim-5})) and this will affect the analysis.

 It was shown in \cite{O-M1}   that  for most system parameters a single piezoelectric beam is not exactly controllable in the natural energy space if magnetic effects are considered.  This is quite different from the conclusion for models without magnetic effects. The paper \cite{Tucs} studied exact controllability of beam-patch system without magnetic effects, and establishes the space of exact controllability strongly depending on the location of the patches. This paper and \cite{O-M1} suggest that the exact controllability of the beam-patch system depends on not only  the location of the patches but also the system parameters. This question is being studied.

%For the future, we aim to extend these beam models to plate models by using Von-Karman  and Reissner-Mindlin approaches.
%%%%%%%%%%%%%%%%%%%%%%%%%%%%%%%%%%%%%%%%%%%%%%%%%%%%%%%%%%%%%%%%%%%%%%%%%%%%%%%%


\begin{thebibliography}{99}

\bibitem{Banks-Smith} H.T. Banks, R.C. Smith, Y. Wang, {\sl Smart material structures: Modelling, Estimation and Control}, Mason, Paris; 1996.
\bibitem{Dietl} J.M. Dietl, A.M. Wickenheiser, E. Garcia, {A Timoshenko beam model for cantilevered piezoelectric energy harvesters,} {\sl{Smart Mater. Struct.}}, vol. 19 (055018), 2010.
\bibitem{E-Inman} A. Erturk, D. Inman, A distrbiuted parameter model for cantilever model for piezoelectric energy harvesting fom base excitations, {\sl{J. Vib. Acoust.,}} vol. 130 (041002), 2008.
\bibitem{Hansen} S. W. Hansen, Analysis of a Plate with a Localized Piezoelectric Patch, {\sl{Proceeedings of the $37^{\text{th}}$ Conference on Decision and Control,}}  Tampa, Florida, 1998, pp. 2952--2957.
\bibitem{K-M-M2} B. Kapitonov, B. Miara, and G. P. Menzala, Stabilization of a layered 3--D body by boundary dissipation, {\sl{ESAIM:COCV}}, vol. 12, 2006, pp. 198--215.

  \bibitem{Lagnese-Lions} J.E. Lagnese, J.-L. Lions, {\sl Modeling Analysis and Control of Thin Plates,} Masson, Paris; (1988).

\bibitem{L-M} I. Lasiecka and B. Miara, Exact controllability of a 3D piezoelectric body, {{\sl C. R. Math. Acad. Sci. Paris}}, vol. 347, 2009, pp. 167--172.
\bibitem{Lee} P.C.Y. Lee,  A variational principle for the equations of piezoelectromagnetism in elastic dielectric crystals,{\sl{ Journal of Applied Physics,}} vol 69 (11), (1991), pp. 7470--7473.	

\bibitem{O-M} K.A. Morris and  A.\"{O}. \"{O}zer, Strong stabilization of piezoelectric beams with magnetic effects, {\sl{Proceeedings of the $52^{\text{nd}}$ Conference on Decision and Control,}}, Florence, Italy, 2013, pp. 3014--3019.

\bibitem{O-M1} K.A. Morris and  A.\"{O}. \"{O}zer, Modeling and stabilizability of voltage-actuated piezoelectric beams with magnetic effects, under revision in  SIAM J. Cont. Optim.

%\bibitem{Scott} W. T. Scott, {\sl{Approximation to real irrationals by certain classes of rational fractions,}} Bull. Amer. Math. Soc. vol. 46 (1940) pp. 124-129.

\bibitem{Ronkanen} P. Ronkanen, P. Kallio, M. Vilkko, H.N. Koivo,  Displacement Control of Piezoelectric Actuators Using Current and Voltage, {\sl IEEE/ASME Trans. Mechatronics}, vol 16 (1), 2011, pp. 160-166.
\bibitem{Rogacheva} N. Rogacheva, {\sl{ The Theory of Piezoelectric Shells and Plates}}, Boca Raton, FL: CRC Press; 1994.
%\bibitem{Sene} A. Sen\`{e}, {\sl{Modelling of piezoelectric static thin plates   ,}} { Asymptotic Analysis},
%(25-1) (2001), pp. 1--20.
\bibitem{Smith} R.C. Smith, {\sl Smart Material Systems}, Society for
Industrial and Applied Mathematics; 2005.
\bibitem{Tiersten} H.F. Tiersten, {\sl {Linear Piezoelectric Plate Vibrations} }, Plenum Press, New York;1969.
%\bibitem{Trigg0} \newblock R. Triggiani, {\newblock\emph{On the stabilizability  problem in Banach space,}} \newblock {J. Math. Anal. Appl.} (3), \textbf{52}  (1975), 383--403.
\bibitem{Tucs} M. Tucsnak, Regularity and exact controllability for a beam with piezoelectric
actuators, {\sl{ SIAM J. Cont. Optim.,}} vol. 34, 1996, pp. 922--930.
%\bibitem{Weiss-Tucsnak} M. Tucsnak and G. Weiss, {\sl{Observation and Control for Operator Semigroups}}, Birkhäuser Verlag, Basel (2009).

%\bibitem{Weiss} M. Tucsnak, G. Weiss,  {\sl{How to get a conservative well-posed linear system out of thin air, Part II: controllability and stability,}} SIAM J. Cont. Optim. (42-3) (2003), pp. 907-935.

\bibitem{Tzou} H.S. Tzou, {\sl{Piezoelectric shells, Solid Mechanics and Its applications 19,}} Kluwer Academic, The Netherlands; 1993.
%\bibitem{Wang-Guo} J-M Wang, B-Z Guo {\sl On the stability of swelling porous elastic soils with fluid
%saturation by one internal damping,} IMA Journal of Applied Mathematics (71)(2006), pp. 565-582.
%\bibitem{Young} R.M. Young, {\sl{An Introduction to Nonharmonic Fourier Series,}} New York, Academic Press, (1980).
    \bibitem{Yang} J. Yang, {\sl{An Introduction to the Theory of Piezoelectricity,}} Springer, New York; 2005.
    \bibitem{Yang1} J. Yang, A review of a few topics in piezoelectricity, {\sl{ Appl. Mech. Rev.,}} vol. 59, 2006, pp. 335–-345.
\bibitem{Zhang} C-G Zhang, Regularity and exact controllability for the Timoshenko beam with Piezoelectric actuator,  Rocky Mountain J. Math., vol. 41 (3), 2011, pp. 999--1010.
\end{thebibliography}
\end{document}